\documentclass[10pt,reqno]{amsart}
     \makeatletter
     \def\section{\@startsection{section}{1}%
     \z@{.7\linespacing\@plus\linespacing}{.5\linespacing}%
     {\bfseries
     \centering
     }}
     \def\@secnumfont{\bfseries}
     \makeatother
\newtheorem{theorem}{Theorem}
\newtheorem{lemma}{Lemma}

\theoremstyle{definition}

\theoremstyle{remark}
\newtheorem{remark}{Remark}
\numberwithin{equation}{section}
\usepackage{natbib}
\setcitestyle{authoryear,open={(},close={)}} 
\usepackage{xcolor}
\usepackage[normalem]{ulem} 
\usepackage{amsmath}
\usepackage{amsthm}
\usepackage{amsfonts}
\usepackage{amssymb}
\usepackage{mathrsfs}
\usepackage{subfig,graphicx}
\usepackage{tikz-cd}
\usepackage{lmodern}

\begin{document}

\title[On Minimal Predictable Intensity]{On Minimal Predictable Intensity of Point Processes}

\author{Haoming Wang}
\address{Haoming Wang: School of Mathematics, Sun Yat-sen University, Xingang Road No. 135, Guangzhou, 510275, China}
\email{wanghm37@mail2.sysu.edu.cn}
\urladdr{https://blueairm.github.io}



\begin{abstract}
An adapted, right‐continuous, non-decreasing, integer‐valued process with unit jumps and starting at zero has a minimal predictable intensity if and only if it is a standard Poisson process under an absolutely continuous transformation of measures.
\end{abstract}

\maketitle

\section{Introduction}  

This paper studies a problem arising in the change of time on semimartingales by transformation of measures. For the continuous semimartingale, \cite{monroe1978aop} proved that any continuous local martingale is a time-changed Brownian motion, although possibly under an enlargement of the filtration. By the Zorn lemma, a minimal filtration always exists in such a problem. For the Brownian motion, this problem is resolved by many authors, such as \cite{kakutani1948equivalence}, \cite{girsanov1960}, \cite{clark1970representation}, \cite{novikov1972identity}, and \cite{kazamaki1994continuous}. In parallel with the Brownian motion, the Poisson process is another prototype in many stochastic applications. In fact, the point process, which is an adapted, right‐continuous, non-decreasing, integer‐valued stochastic process with unit jumps and starting at zero, can be represented as a time-changed Poisson process. However, an enlargement of the filtration is still required. \cite{bremaud1972} once asked when the minimal filtration is natural and conjectured that if and only if it is a standard Poisson process under an absolutely continuous transformation of measures. In this article, we answer this question in the affirmative sense.

Although on many occasions, an enlargement of filtration helps us represent point processes, it brings some troubles in the analysis of such processes, e.g., including more null sets. Converse to the filtration expansion problem, like expositions in \cite{protter2005stochastic}, \cite{jeulin2006semi}, and \cite{DITELLA2021103}, the minimal filtration was first studied by \cite{jacod1975}, and subsequently by \cite{bremaud1981point}, \cite{daley2003introduction} and \cite{karr2017point}. Their methods are all based on the hazard function of conditional expectation. For example, Theorem 14.1 in \cite{daley2003introduction} and Theorem 2.31 in \cite{karr2017point}. In this method, it is unavoidable to assume the null complete assumption in the sense of \cite{dellacherie1972capacites}, i.e., all terminal null sets are contained in the initial filter. However, the explicit construction based on null completeness limits the practicality of the minimal filtration. In our paper, we adopt a function-analytic approach to establish the existence of a minimal predictable intensity using the Hahn-Banach theorem, which only involves the property of closed subspaces of semimartingales. This method helps us get rid of unnecessary null completeness, so that we can obtain a natural restriction of filtration.

Further applications to the Cox and Hawkes processes implemented the main theorem. 
Based on our main Theorems and the Lebesgue decomposition theorem, we establish a canonical decomposition Theorem \ref{thm: decomposition theorem} for the Cox process. With application to the Hawkes process, Theorem \ref{thm: Hawkes Poisson process} fills existing gaps between Poisson processes and self-exciting processes.

\section{Notations and Main Results}

The present work studies a class of stochastic processes defined in the real half line $[0,\infty)$. Let $\mathscr{F}_{t}$ be a monotone family of $\sigma$-algebras of subsets of an abstract sample space $\boldsymbol{\varOmega}$, $\mathscr{F}_{s} \subseteq \mathscr{F}_{t} \subseteq \mathscr{F}$ ($t > s$), bounded by the $\sigma$-algebra $\mathscr{F}$ of all events of possible outcome. Let us denote $\{\mathscr{F}_t: t\geq 0\}$ as $\boldsymbol{F}$. We are only interested in the case of a process
\[\underline{\underline{N}} = \{N(\omega,t), \boldsymbol{F},  \boldsymbol{P}\}\]
such that $N(\omega,t)$ is a piecewise constant function on $[0,\infty)$ for all $\omega$, where $\boldsymbol{P}$ is a probability measure defined on the $\sigma$-algebra of events $\mathscr{F}$. Here it is assumed that the set $\{\omega: N(\omega, t) \in \varGamma\} \in \mathscr{F}_{t}$ for any $t \geq 0$ and any subset $\varGamma$ on the real line.

The process $\underline{\underline{N}}$ is called Poisson of rate $\lambda$ if it is 
\begin{enumerate}
    \item[(P0)] (adapted) $\{\omega: N(\omega, t) \in \varGamma\} \in \mathscr{F}_{t}$ for any $t \geq 0$ and any Borel set $\varGamma$ on the real line;
    
    \item[(P1)] (right-continuous) $\lim\limits_{t\to s+}N(\omega,s) = N(\omega,s)$ ($f(x+)$ is used for the right limit $\lim\limits_{y\to x+}f(y)$ and $f(x-)$ for the left limit $\lim\limits_{y\to x-}f(y)$); 

    \item[(P2)] (non-decreasing) $N(\omega, t) \geq N(\omega, s), t > s$; 

    \item[(P3)] (integer-valued) $N(\omega,t)\in \{0,\pm 1,\pm 2,\dots\}$; 
    
    \item[(P4)] (starting from zero) $N(\omega,0)=0$;

    \item[(P5)] (with unit jumps) $\Delta N(\omega,s) = N(\omega,s) - N(\omega,s-)  \in \{0,1\}$,\\
    for fixed $\omega$; and if it has
    \item[(P6)] (Poisson law) $\int_{\varOmega} \exp [iu (N(\omega,t) - N(\omega,s))] \boldsymbol{P}(d\omega) = \exp[\lambda (t-s)(e^{iu}-1)]$;
    \item[(P7)] (independent increment) $N_{t} - N_{s} \perp \mathscr{F}_{s}$ for $0\leq s \leq t$.
\end{enumerate}

In fact, from a theorem of \cite{renyi1967remarks}, the pairwise independence (this symbol $\perp$ is used for this.) of the Poisson process implies the $n$-component independence. Define $T_{n}(\omega) = \inf \{t\geq 0: N(\omega,t) \geq n \}$. It is well-known that a Poisson process can be represented as 
\begin{equation}
    N(\omega,t) = \sum_{n=1}^{\infty} \mathbf{1}\{T_{n}(\omega) \leq t\},
    \label{eq: Poisson process}
\end{equation}
such that $U_{n+1}(\omega) = T_{n + 1}(\omega) - T_{n}(\omega)$ where $n = 1,2,\dots$, follow the independent identically distributed exponential distribution with mean $1/\lambda$.

The definition of \eqref{eq: Poisson process} can be generalised to any such sequence of stopping times $0 ( = T_{0}(\omega) ) \leq T_{1}(\omega) \leq T_{2}(\omega) \leq \dots$ with $T_{n}(\omega)$ tending to infinity and $\{\omega: T_{n}(\omega) \leq t\} \in \mathscr{F}_{t}$ for all $t$ and $n$. It is adapted since it is a sum of indicators of $\{\omega: T_{n}(\omega) \leq t\} \in \mathscr{F}_{t}$ for fixed $t$; right-continuous since the indicators are non-negative and right-continuous for fixed $\omega$. It is obviously non-decreasing and integer-valued. However, it does not necessarily start from zero or have unit jumps either since two such $T_{n} (\omega)$ and $T_{n'}(\omega), n \neq n'$ may be equal on some subset $\boldsymbol{\varOmega}^{\prime} \in \mathscr{F}$ of $\boldsymbol{\varOmega}$ such that $\boldsymbol{P}(\boldsymbol{\varOmega}^{\prime})>0$. Except for this situation, the process \(\underline{\underline{N}}\) satisfying (P0)-(P5) will be called a simple point process. We say that an adapted process is predictable if it is jointly measurable in both variables with respect to the smallest $\sigma$-algebra on $\boldsymbol{\varOmega}\times [0,\infty)$ generated by all left-continuous adapted processes. The compensator of a simple point process $\underline{\underline{N}}$ is the unique predictable non-decreasing process  \( A(\omega,t) \) with \( A(\omega, 0) = 0 \) such that \( N(\omega,t) - A(\omega,t)\) is an $\boldsymbol F$-local martingale. A simple point process is locally of integrable variation (See this definition on Page 66 before Theorem 22 of Chapter III in \cite{protter2005stochastic}.), so its compensator always exists. For example, the compensator of a Poisson process is \( A(\omega,t) = \lambda t \). 

Coming to a simple point process $\underline{\underline{N}} = \{N(\omega,t),\boldsymbol{F}, \boldsymbol{P}\}$, if \(\boldsymbol{F}\) coincides the natural filtration generated by all $N(\cdot,s), s\leq t$, we will say that simple point process $\underline{\underline{N}} = \{N(\omega,t),\boldsymbol{F}, \boldsymbol{P}\}$ is non-anticipative and that its compensator \( A(\omega,t)\) is natural. By the Zorn lemma, there exists a minimal filtration \(\boldsymbol{F} \) in the partially ordered set $(S,\leq)$ of all filtration $\boldsymbol F$ satisfying
\begin{enumerate}
    \item[(S1)] $N \text{ is }\boldsymbol{F}\text{-adapted},$
    \item[(S2)] $A \text{ is }\boldsymbol{F}\text{-predictable},$ 
    \item[(S3)] $N - A \text{ is an }\boldsymbol{F}\text{-local martingale}$,
\end{enumerate}
where $\boldsymbol{G} \leq \boldsymbol{F} \Leftrightarrow \mathscr{G}_{t}\in\boldsymbol{G} \subseteq \mathscr{F}_{t} \in \boldsymbol{F},$ for all $t \geq 0$. This is because the non-empty intersection of a descending chain of filtration in $S$ again satisfies (S1)-(S3), thus yielding a lower bound, so the Zorn lemma applies.  \cite{bremaud1972} asked a question: Should the minimal filtration $\boldsymbol F \neq \boldsymbol F^N$ if there exists a probability measure $\boldsymbol{P}_0$ relative to which the point process $\underline{\underline{N}} = \{N(\omega,t),\boldsymbol{F}, \boldsymbol{P}\}$ is absolutely continuous and $\underline{\underline{N}}^0 = \{N(\omega,t),\boldsymbol{F}, \boldsymbol{P}_0\}$ is the standard Poisson process? Here we try to give a brief answer to his question that it is impossible. 

Recall that a simple point process \( \underline{\underline{N}} \) is not necessarily orderly, by which we mean that for any stopping time \( \tau(\omega) \geq 0 \),  
\begin{equation}\label{eq: orderly}  
\lim_{h \to 0} \frac{1}{h} \boldsymbol{P}(N(\cdot,\tau(\cdot)+h) - N(\cdot,\tau(\cdot)) > 1 \mid \mathscr{F}_{\tau})(\omega) = 0,
\end{equation}  
while orderliness implies simplicity (e.g., the mixed Poisson process of Exercise 3.3.2 in \cite{daley2003introduction}). Here, $\mathscr{F}_{\tau}$ denotes the stopped $\sigma$-subalgebra of $\mathscr{F}$,
\[\mathscr{F}_{\tau} = \{A: A \cap \{\omega: \tau(\omega) \leq t\} \in \mathscr{F}_t,  \text{ for all } t \geq 0\},\] 
The intensity \( \lambda(\omega,t) \)  of an orderly point process $\underline{\underline{N}}$ is a non-negative process satisfying for any predictable stopping time \( \tau(\omega) \geq 0 \),
\begin{equation}\label{eq: intensity}  
\lim_{h \to 0} \frac{1}{h} \boldsymbol{P}(N(\cdot,\tau(\cdot)+h) - N(\cdot,\tau(\cdot)) = 1 \mid \mathscr{F}_{\tau}) (\omega) = \lambda(\omega,\tau(\omega)).
\end{equation}  
Here, we use a different definition from \cite{bremaud1972}, who evaluates \eqref{eq: orderly}  and \eqref{eq: intensity} for fixed $t \ge 0$. The advantage of our definition lies in that it requires a stronger uniqueness of indistinguishability for two such processes $\lambda(\omega,t)$ and $\lambda'(\omega,t)$, i.e.,
\[\boldsymbol{P}(\lambda(\cdot,t) \neq \lambda'(\cdot,t) \text{ for some } t \ge 0) = 0.\]
This is guaranteed by the predictable section theorem (Theorem 85 of Chapter IV  in \cite{Dellacherie1980probabilities}), while in Br\'emaud's definition, only version uniqueness is assumed, i.e.,
\[\boldsymbol{P}(\lambda(\cdot,t) \neq \lambda'(\cdot,t)) = 0,  \text{ for all } t \ge 0.\]

Let $\boldsymbol{F}$ be an arbitrary filtration and $\boldsymbol{F}^N$ the natural filtration generated by a simple point process $\underline{\underline{N}}$. Also, let $\mathscr{F}_{\infty} 
$ be the smallest $\sigma$-subalgebra of $\mathscr F$ containing all $\mathscr{F}_t$, $t \geq 0$.  We shall assume $\lambda(\omega,t)$ satisfies
    \begin{itemize}
        \item[(A0)] $\lambda(\omega,t)$ is $\boldsymbol{F}$-adapted;
        \item[(A1)] $\lambda(\omega,t)$ is measurable for both $\omega$ and $t$;
        \item[(A2)] $\lambda(\omega,t)$ is strictly positive for fixed $\omega$ and $t$;        
        \item[(A3)] $\int_{0}^{t} \lambda(\omega,s) ds < \infty$ almost surely for each $t$;
        \item[(A4)] $\int_{0}^{t} \log \lambda(\omega,s) N(\omega, ds)$ exists almost surely for each $t$.
    \end{itemize}
    From $\lambda(\omega,t)$ we define $M(\lambda)$, $A(\lambda)$, and $\zeta_{s}^{t}(\lambda)$:
    \begin{equation*} 
    \begin{aligned}
        M(\lambda) (\omega,t) = N(\omega,t) - A(\lambda)(\omega), \text{ where } A(\lambda)(\omega,t) = \int_{0}^{t} {\lambda(\omega,s)} ds,\\
        \text{ and } \zeta_{s}^{t}(\lambda)(\omega) = \int_{s}^{t} [1 - \lambda(\omega,u)]du + \int_{s}^{t} \log \lambda(\omega,u) N(\omega,du).
    \end{aligned}
    \end{equation*}
    For $t \in [0,\infty]$, let us set 
        \begin{equation}
            \boldsymbol{Q}(d\omega) = \exp[\zeta_0^t(\lambda)(\omega)]\boldsymbol{P}(d \omega), \text{ on } \mathscr{F}_{t}.
            \label{eq: girsanov transformation}
        \end{equation}
    \begin{theorem}    \label{thm: main theorem 1}
    Suppose $\underline{\underline{N}} = \{N(\omega,t),\boldsymbol{F}, \boldsymbol{P}\}$ is a Poisson process of rate 1 and $\lambda(\omega,t)$ satisfies \textnormal{(A0)-(A4)}. Take $t = \infty$ in \eqref{eq: girsanov transformation}. 
    We assume that $\boldsymbol{Q}(\boldsymbol{\varOmega}) = 1.$ For any $\boldsymbol{G}$ such that $\boldsymbol{F}^{N} \leq \boldsymbol{G} \leq \boldsymbol{F}$, there exists a unique $\boldsymbol{G}$-predictable process $\lambda^{\boldsymbol{G}}(\omega,t)$ up to indistinguishability such that $\underline{\underline{M(\lambda^{\boldsymbol{G}})}} = \{M(\lambda^{\boldsymbol{G}})(\omega,t),\boldsymbol{G}, \boldsymbol{Q}\}$ is a martingale. In particular, if $\lambda(\omega,t)$ is $\boldsymbol{G}$-predictable, then it is indistinguishable to $\lambda^{\boldsymbol{G}}(\omega,t)$.
    \end{theorem}  
        
\begin{remark}
    From the section theorem and stopping theorem (Theorem 43 of Chapter VI in \cite{Dellacherie1980probabilities}), any non-negative $\boldsymbol{F}$-measurable process $\lambda(\omega,t)$ has a $\boldsymbol{G}$-predictable projection ${}^p\lambda(\omega,t)$ such that for any finite non-negative predictable stopping time $\tau(\omega)\ge 0$, 
        \begin{equation}
            \boldsymbol{E_P}[\lambda(\omega,\tau(\omega))\mid \mathscr{G}_{\tau-}](\omega) =  {}^p\lambda(\omega,\tau(\omega))
        \end{equation}
        where $\mathscr{G}_{0-} = \mathscr{G}_{0}$ and  
        \[\mathscr{G}_{\tau-}=\sigma\left(\left\{ A\cap\{\omega:t<\tau(\omega)\}\colon t\ge 0,A\in\mathscr{G}_t \right\}\cup\mathscr{G}_0\right). \]
        We can easily check $\mathscr{G}_{\tau-}$ is the smallest $\sigma$-algebra containing all $\mathscr{G}_{\tau_{n}}$ for any sequence of non-negative non-decreasing stopping times $\tau_n(\omega)$ tending to $\tau(\omega)$ (Theorem 6 of Chapter III in \cite{protter2005stochastic}). Thus, Theorem \ref{thm: main theorem 1} can be restated as follows: ${}^p\lambda(\omega,t)$ is the unique $\boldsymbol{G}$-predictable process up to indistinguishability such that $\underline{\underline{M({}^p\lambda)}} = \{M({}^p\lambda)(\omega,t),\boldsymbol{G}, \boldsymbol{Q}\}$ is a martingale.
\end{remark}

\begin{remark}
    The converse of Theorem \ref{thm: main theorem 1} is also true, which will be stated in Theorem \ref{thm: main theorem 2} of Section 4. 
\end{remark}

Hereafter, we will use these notions for short: $N_{t}(\omega) = N(\omega,t)$ sometimes with $\omega$ omitted, and $\boldsymbol{E_P}[Z] = \int_{\Omega} Z(\omega)\boldsymbol{P}(d\omega)$ sometimes with the subscript $\boldsymbol{P}$ left out when the underlying probability measure $\boldsymbol{P}$ remains clear. 

\section{Proof of  Main Theorem \ref{thm: main theorem 1}}

        
\begin{lemma} \label{thm: watanabe-bremaud} Suppose $\underline{\underline{N}} = \{N(\omega,t),\boldsymbol{F}, \boldsymbol{P}\}$ satisfies \textnormal{(P0)-(P5)} and
\begin{equation}
    \int_{\boldsymbol{\varOmega}}N(\omega,t)\boldsymbol{P}(d\omega) < \infty, \text{ for each } t \geq 0,\label{eq: integrable process}
\end{equation}
Suppose also $\underline{\underline{\lambda}} = \{\lambda(\omega,t),\boldsymbol{F}, \boldsymbol{P}\}$ satisfies \textnormal{(A0)-(A4)}. The following statements are equivalent. 
	\begin{enumerate}
	    \item For any non-negative $\boldsymbol{F}$-predictable process $C(\omega,t)$,
		$$\int_{\boldsymbol{\varOmega}}\boldsymbol{P}(d\omega)\left[\int_{0}^{\infty}C(\omega,t)N(\omega,dt)\right] = \int_{\boldsymbol{\varOmega}}\boldsymbol{P}(d\omega)\left[\int_{0}^{\infty}C(\omega,t)A(\omega,dt)\right].$$	
            \item $N(\omega,t) - A(\omega,t)$ is an $\boldsymbol{F}$-martingale under $\boldsymbol{P}$.
            \item $N(\omega,t) - A(\omega,t)$ is an $\boldsymbol{F}$-local martingale under $\boldsymbol{P}$.
	\end{enumerate}
	\end{lemma} 
\begin{proof}[Proof of Lemma \ref{thm: watanabe-bremaud}]
(3) $\Rightarrow$ (1). If $N(\omega,t) - A(\omega,t)$ is a local martingale, then $A(\omega,t)$ is integrable since in \eqref{eq: integrable process} the $N(\omega,t)$ is required to be integrable for each $t \geq 0$. Let's choose a sequence of non-negative stopping time $\tau_{n}(\omega)$ tending to infinity and a predictable generator $$C(\omega,u) = \mathbf{1}\{\omega\in \varGamma\}\mathbf{1}\{s < u \leq t\}$$ for some $\varGamma \in \mathscr{F}_{s}$, 
        \begin{equation}
		  \begin{aligned}
                \boldsymbol{E}[(N_{t\land \tau_{n}} & - N_{s\land \tau_{n}});\varGamma] = \boldsymbol{E}[\int_{0}^{\infty} C_{u} dN_{u\land \tau_{n}}] \\
			& = \boldsymbol{E}[\int_{0}^{\infty} C_{u} dA_{u\land \tau_{n}}] = \boldsymbol{E}[(A_{t\land \tau_{n}} - A_{s\land \tau_{n}});\varGamma],
		  \end{aligned}
			\label{eq: martingale}
	\end{equation}
        where the first and the last expectation in \eqref{eq: martingale} is taken over a subset $\varGamma$ of $\boldsymbol{\varOmega}$. Letting $n$ tend to $\infty$, (1) is thereby obtained by the monotone class theorem. If (\ref{eq: martingale}) holds for any predictable process $C(\omega,u)$, it proves that $N(\omega,t)$ is a martingale again by the integrability condition \eqref{eq: integrable process} of $N(\omega,t)$. The converse (1) $\Rightarrow$ (2) is also true. Since any martingale is a local martingale, (2) $\Rightarrow$ (3) holds.
\end{proof}

            \begin{lemma}\label{lem: exp martingale} Assume the same conditions of Theorem \ref{thm: main theorem 1}. If we put $Z(\omega,t) = \exp[\zeta_{0}^{t}(\omega)]$, 
        then $\boldsymbol{E}[Z_{\infty}] =1$ implies $\underline{\underline{Z}} = \{Z(\omega,t),\boldsymbol{F}, \boldsymbol{P}\}$ is a uniformly integrable martingale, i.e., $\underline{\underline{Z}}$ is a martingale and 
        \begin{enumerate}
            \item $\exists C>0, \sup_{t}\boldsymbol{E}[|Z_{t}|]\leq C$;\\
            \item $\forall \varepsilon >0, \exists \delta >0$ such that $\forall \varGamma \in \mathscr{F}_{\infty}$, $P(\varGamma)\leq \delta$, $\sup_t\boldsymbol{E} [|Z_t|; \varGamma]< \varepsilon$.
        \end{enumerate}
        \end{lemma}
        \begin{proof}[Proof of Lemma \ref{lem: exp martingale}]
        Let us rewrite $Z(\omega,t) = \exp[\zeta_{0}^{t}(\omega)]$ as 
            \[\begin{aligned}
                 Z_t = \exp \left[\int_{0}^{t} \log \lambda(\omega,s) N(\omega,ds)\right] \Big/ \exp\left[\int_{0}^{t} (\lambda(\omega,t)-1)ds \right] (= X_t/Y_t). 
            \end{aligned}\]
                We are going to show that $\underline{\underline{Z}}$ is a martingale. By taking the differential,
                we find from $\underline{\underline{N}}$ is simple that 
                \[
                \begin{aligned}
                    dZ = YdX + XdY = Y(dX - X(\lambda -1)ds) (=YdM),\\
                M(\omega,t) = X(\omega,0) + \int_{0}^{t} [X(\omega,s) - X(\omega,s-)][N(\omega,ds) - ds].
                \end{aligned}\]
                Here, $\underline{\underline{M}}$ is a martingale. $\underline{\underline{Z}}$ is a local martingale since $\underline{\underline{Y}}$ is predictable. Let its associated stopping time be $\tau_{n} (\omega) \geq 0$. By the Fatou lemma,
                \[\boldsymbol{E} [Z_{t}\mid\mathscr{F}_{s}] \leq \liminf_{n\to \infty} \boldsymbol{E}[Z_{t\land \tau_{n}}\mid\mathscr{F}_{s}] = \liminf_{n\to \infty} Z_{s\land \tau_{n}} = Z_{s}, \quad t \geq s\]
                Hence, $\underline{\underline{Z}}$ is a non-negative super-martingale. 
                As $\boldsymbol{E}[Z_{t}|\mathscr{F}_{s}]$ takes values between $0$ and $Z_{s}$, the common expectation $\boldsymbol{E} [Z_{\infty}] = \boldsymbol{E} [Z_{0}] = 1$ implies the equality holds. Thus, $\underline{\underline{Z}}$ is a martingale. 
            \end{proof}
        The uniform integrability follows from the following basic lemma, the proof of which can be found in Theorem 11 of Chapter I in \cite{protter2005stochastic}.
        \begin{lemma}\label{lem: revuz yor}
           Suppose $\underline{\underline{Z}} = \{Z(\omega,t),\boldsymbol{F}, \boldsymbol{P}\}$ is a martingale. The following three conditions are equivalent.
            \begin{enumerate}
                \item There exists a random variable $Z_{\infty}$ such that $\lim\limits_{t\to\infty}\boldsymbol{E}|Z_t - Z_{\infty}|= 0$.
                \item There exists a random variable $Z_{\infty}$ such that $\boldsymbol{E} [Z_{\infty}]<\infty$ and $Z(\omega,t) = \boldsymbol{E}[Z_{\infty}|\mathscr{F}_{t}](\omega)$.
                \item $\underline{\underline{Z}} = \{Z(\omega,t),\boldsymbol{F}, \boldsymbol{P}\}$ is uniformlly integrable.
            \end{enumerate}
        \end{lemma}

        \begin{lemma}  Under \textnormal{(A0)-(A4)},  $\underline{\underline{M(\lambda)}} = \{M(\lambda)(\omega,t),\boldsymbol{F}, \boldsymbol{Q}\}$ is a martingale.   \label{lem: compensated martingale}        \end{lemma}
        \begin{proof}[Proof of Lemma \ref{lem: compensated martingale}]
            By Lemma \ref{thm: watanabe-bremaud}, we need to show that for any non-negative $\boldsymbol{F}$-predictable process $C(\omega,t)$, 
        \begin{equation}\begin{aligned}
				\boldsymbol{E}_{\boldsymbol{Q}}[\int_{0}^{\infty} C_t dN_t] =  \boldsymbol{E}_{\boldsymbol{Q}}[\int_{0}^{\infty}C_t \lambda_td{t}].
        			\end{aligned}
                    \label{eq: girsanov}
		\end{equation}
        By the definition of \(\boldsymbol{Q}\), 
        \eqref{eq: girsanov} is equivalent to
        \begin{equation}\begin{aligned}
			\boldsymbol{E}_{\boldsymbol{P}}[Z_\infty\int_{0}^{\infty} C_t dN_t] =  \boldsymbol{E}_{\boldsymbol{P}}[Z_\infty\int_{0}^{\infty}C_t \lambda_td{t}],
        			\end{aligned}
                    \label{eq: girsanov 2}
	\end{equation}
        Combining the Fubini theorem and Lemma \ref{lem: exp martingale}, it suffices to show
        \begin{equation}\begin{aligned}
				\boldsymbol{E}_{\boldsymbol{P}}[\int_{0}^{\infty} Z_tC_t dN_t] =  \boldsymbol{E}_{\boldsymbol{P}}[\int_{0}^{\infty}Z_tC_t \lambda_td{t}].
        			\end{aligned}
                    \label{eq: girsanov 3}
	\end{equation}
        By the construction of \(Z(\omega,t)\), $Z(\omega,t) = Z(\omega,t-)  \lambda(\omega,t)$ at jump times $t = T_{n}(\omega),$ and $Z(\omega,t) = Z(\omega,t-)$ otherwise.
        The LHS of \eqref{eq: girsanov 3} then reduces to
        \begin{equation}\begin{aligned}
				\boldsymbol{E}_{\boldsymbol{P}}[\int_{0}^{\infty} Z_{t-}C_t \lambda_tdN_t] =  \boldsymbol{E}_{\boldsymbol{P}}[\int_{0}^{\infty}Z_{t-}C_t \lambda_td{t}]
        			\end{aligned}
                    \label{eq: girsanov 4}
	\end{equation}
        Here, we use the fact that \(N(\omega,t) - t\) is a $\boldsymbol{F}$-martingale under $\boldsymbol{P}$. \eqref{eq: girsanov 4} yields the desired result \eqref{eq: girsanov 3} because any countable set on the real half-line has Lebesgue measure zero. 
    \end{proof}

         \begin{lemma}
            \label{lem: reduction}
            If Theorem \ref{thm: main theorem 1} holds for $\boldsymbol{G}$, then it holds for any $\boldsymbol{G}^{\prime}$ satisfying $\boldsymbol{G}\leq \boldsymbol{G}^{\prime} \leq \boldsymbol{F}$.
        \end{lemma}
        \begin{proof}[Proof of Lemma \ref{lem: reduction}] Direct from Lemma \ref{lem: revuz yor}.
        \end{proof}

        Therefore, we could focus on proving Theorem \ref{thm: main theorem 1} only for $\boldsymbol{G} = \boldsymbol{F}^N$.
        The following lemma is a version of Theorem 1.21 in Chapter III of \cite{jacod2013limit}. We omit the proof.
        \begin{lemma} \label{lem: jacod} Assume the same conditions in Lemma \ref{thm: watanabe-bremaud}.  Let $\boldsymbol{P}$ and $\boldsymbol{Q}$ be two probability measures on $(\boldsymbol{\varOmega}, \mathscr{F}_{\infty}^{N})$ with respect to which $N(\omega.t)$ has common $\boldsymbol{F}^{N}$-compensator. Then $\boldsymbol{P}$ and $\boldsymbol{Q}$ coincide on $\mathscr{F}_{\infty}^{N}$.		\end{lemma} 
        
        We are going to prove the predictable representation property of point processes via the Hahn-Banach theorem. This proof is essentially due to \cite{2010159}.

        \begin{lemma}\label{thm: mrp} Assume the same conditions in Theorem \ref{thm: main theorem 1} and also that one of the statements in Lemma \ref{thm: watanabe-bremaud} holds. If $\underline{\underline{L}} = \{L(\omega,t),\boldsymbol{F}^{N}, \boldsymbol{P}\}$ is a uniformly integrable martingale, then there exists a unique $\boldsymbol{F}^{N}$-predictable process $H(\omega,t)$ up to indistinguishability such that except a $\boldsymbol P$-null set,
            \begin{equation}\label{eq: uniformly integrable Mt}
                L(\omega,t) = L(\omega,0) + \int_{0}^{t}H(\omega,s)[N(\omega,ds) - A(\omega,ds)], \text{ \, for each $t$}.
            \end{equation}
            \end{lemma}
        \begin{proof}[Proof of Lemma \ref{thm: mrp}] Let $L^1(\boldsymbol P)$ consists of the processes $L(\omega,t)$ in \eqref{eq: uniformly integrable Mt} with $H(\omega,t)$ running over all $\boldsymbol{F}^{N}$-predictable processes. Clear, $L^1(\boldsymbol P) \subseteq M_{\textrm{loc}} (\boldsymbol P)$, the set of all local martingales under $\boldsymbol P$. We aim to prove that the inclusion relation is an equality. If not, since $L^1(\boldsymbol P)$ is a properly closed subspace, by the Hahn-Banach theorem, there exists a continuous linear function $f \in M_{\textrm{loc}} (\boldsymbol P)^*$ such that $f(L) = 0$, for all $L \in L^1(\boldsymbol P)$, not vanishing identically. But from classical Hardy-BMO duality in Theorem 10.21 of \cite{he1992}, there exists a non-zero process $K \in BMO$ such that 
        \[\boldsymbol E \left( [L,K]_{\infty}\right) = 0, \text{ for all } L \in L^1(\boldsymbol P).\]
        Note that $BMO \subseteq M^{\infty}_{\textrm{loc}} (\boldsymbol P)$, the set of all bounded local martingales under $\boldsymbol P$, so there exists a stopping time $\tau$ such that $|K_{t\land \tau}| \leq M$. Take for simplicity $M=1$. A simple calculation shows
        \[\boldsymbol E \left( L_\infty K_{0\land \tau}\right) = \boldsymbol E \left( L_\infty K_{\infty\land \tau}\right) = \boldsymbol E \left( [L,K_{t\land \tau}]_{\infty}\right) = \boldsymbol E \left( [L_{t\land\tau},K]_{\infty}\right) = 0\]
        since $L_{t\land \tau} \in L^1(\boldsymbol P)$. Thus, $K_{0} = K_{0\land \tau}  = 0$ since 
        $\mathscr F^N_0 = \{\boldsymbol\emptyset,\boldsymbol \varOmega\}$ and we can take $L \equiv 1$. Define
        \[Z_t = 1 + \frac{1}{2}K_{t\land \tau}, \quad \frac{1}{2} \leq Z_t \leq \frac{3}{2}.\]
        In particular, we have $\boldsymbol E [Z_t] = \boldsymbol E [Z_0] = 1$, and $\boldsymbol E\left(L_\infty Z_{\infty}\right) = 0$ for $ L\in L^1(\boldsymbol P)$ with $L(\omega,0) =0$.
        This shows that the probability measure $\boldsymbol P'$ defined by $d\boldsymbol P' = Z_{\infty} d\boldsymbol P$ satisfies $\boldsymbol P' \sim \boldsymbol P$. By Lemma \ref{thm: watanabe-bremaud}, $\underline{\underline{M}}' = \{M(\omega,t),\boldsymbol{F}^{N}, \boldsymbol{P}'\}$ is a martingale, where $M(\omega,t) = N(\omega,t) -A(\omega,t)$ . From Lemma \ref{lem: jacod}, $\boldsymbol P' = \boldsymbol P$, so $Z_\infty = 1$, that is, $K_{\tau} =0$. Since $\tau$ can be made arbitrarily large, $K \equiv 0$. Contradiction.
        \end{proof}
        
        \begin{lemma}\label{lem: likelihood}
            $L(\omega,t) = \boldsymbol{E}_{\boldsymbol{P}} [d\boldsymbol{Q}/d\boldsymbol{P}| \mathscr{F}^N_t]$ is an $\boldsymbol{F}^N$-uniformly integrable martingale under $\boldsymbol{P}$.
        \end{lemma}

        \begin{proof}[Proof of Lemma \ref{lem: likelihood}] Direct from Lemma \ref{lem: revuz yor}.
        \end{proof}

        \begin{lemma}  Under \textnormal{(A0)-(A4)},  $\underline{\underline{M(\lambda^N)}} = \{M(\lambda^N)(\omega,t),\boldsymbol{F}^N, \boldsymbol{Q}\}$ is a martingale.   \label{lem: compensated martingale N}        \end{lemma}
        
        \begin{proof}[Proof of Theorem \ref{lem: compensated martingale N}] 
        According to the martingale representation Lemma \ref{thm: mrp}, there exists an $\boldsymbol{F}^{N}$-predictable process $H(\omega,t)$ such that
		\begin{eqnarray}
                L(\omega,0) = 1, \, L(\omega,t) = 1 + \int_{0}^{t} H(\omega,s)(N(\omega,ds) -ds).
			\label{eq: representation}
		\end{eqnarray}
        From \eqref{eq: representation} we have $L(\omega,t) = L(\omega,t-)+ H(\omega,t)$ at jump times $t = T_{n}(\omega)$, and $L(\omega,t) = L(\omega,t-)$ otherwise.
        Let us introduce $\tau(\omega) = \inf \left\{t: L(\omega,t) = 0\right\}$ on the extended real half-line $[0,\infty]$, which is strictly positive with possible values at $\infty$ for fixed $\omega$. Furthermore, $\tau(\omega)$ is also an $\boldsymbol{F}^{N}$-stopping time from the measurable section theorem (\cite{Dellacherie1980probabilities}). It follows that 
        \begin{equation}
            \lambda^{N}(\omega,t) = \left[1 + H(\omega,t)L^{-1}(\omega,t-)\right]\mathbf{1}\{\omega: \tau(\omega) > t\}
        \end{equation}
        is a non-negative $\boldsymbol{F}^{N}$-predictable process, which fulfills \eqref{eq: girsanov}. By Lemma \ref{thm: watanabe-bremaud}, $\underline{\underline{M(\lambda^N)}} = \{M(\lambda^N)(\omega,t),\boldsymbol{F}^N, \boldsymbol{Q}\}$ is thereby a martingale.
        \end{proof}
        
        Uniqueness is again from Lemma \ref{thm: watanabe-bremaud}. Combining Lemma \ref{lem: compensated martingale} and Lemma \ref{lem: compensated martingale N}, we arrive at this conclusion.
        \begin{lemma} Assume the same conditions in Theorem \ref{thm: main theorem 1}. $\lambda^N(\omega,t)$ is the unique $\boldsymbol{F}^N$-predictable process such that Theorem \ref{thm: main theorem 1} holds. \label{lem: compensated martingale on FN}         \end{lemma}

\section{Proof of the Converse of Main Theorem \ref{thm: main theorem 1}}
        
        In fact, we also have the converse of Theorem \ref{thm: main theorem 1}.
        
        \begin{theorem}\label{thm: main theorem 2}
        Suppose \(\underline{\underline{N}} = \{N(\omega,t),\boldsymbol{F}^{N}, \boldsymbol{P}\}\) is a Poisson process with rate 1 and \(\lambda(\omega,t)\) is an \(\boldsymbol{F}^{N}\)-predictable process satisfying \textnormal{(A0)-(A4)}. Let \(\boldsymbol{Q}\) be another probability measure 
        such that \(\underline{\underline{M(\lambda)}} = \{M(\lambda),\boldsymbol{F}^{N}, \boldsymbol{Q}\}\) is a martingale. Then \eqref{eq: girsanov transformation} holds for $\lambda(\omega,t)$ for all $t \in [0,\infty]$. \end{theorem}


    
    \begin{proof}[Proof of Theorem \ref{thm: main theorem 2}] According to Lemma \ref{lem: jacod}, it suffices to show that for each $t \ge 0$, the restriction $\boldsymbol{Q}_t$ of $\boldsymbol{Q}$ on ${\mathscr{F}_t}$ coincides with the probability measure $\tilde{\boldsymbol{Q}}_t$ defined by $\tilde{\boldsymbol{Q}}_t(d\omega) = \exp[\zeta_{0}^{t}(\lambda)(\omega)]\boldsymbol{P}(d\omega)$ on $\mathscr{F}^{N}_{t}$. By Lemma \ref{thm: watanabe-bremaud}, it thereby suffices to prove that for any non-negative $\boldsymbol{F}^{N}$-predictable process $C(\omega,s)$ vanishing outside $[0,t]$, 
            \begin{equation}\begin{aligned}
				\boldsymbol{E}_{\tilde{\boldsymbol{Q}}_t}[\int_{0}^{t} C_s dN_s] =  \boldsymbol{E}_{\tilde{\boldsymbol{Q}}_t}[\int_{0}^{t}C_s \lambda_sd{s}].
        			\end{aligned}
                    \label{eq: girsanov t }
		\end{equation} 
        By the definition of \(\tilde{\boldsymbol{Q}}_t\),
        \eqref{eq: girsanov t } becomes
        \begin{equation}\begin{aligned}
			\boldsymbol{E}_{\boldsymbol{P}}[Z_t\int_{0}^{t} C_s dN_s] =  \boldsymbol{E}_{\boldsymbol{P}}[Z_t\int_{0}^{t}C_s \lambda_sd{s}],
        			\end{aligned}
                    \label{eq: girsanov t 2}	\end{equation}
        Combining the Fubini theorem and Lemma \ref{lem: exp martingale}, it suffices to show
        \begin{equation}\begin{aligned}
				\boldsymbol{E}_{\boldsymbol{P}}[\int_{0}^{t} Z_sC_s dN_s] =  \boldsymbol{E}_{\boldsymbol{P}}[\int_{0}^{t}Z_sC_s \lambda_sd{s}].
        			\end{aligned}                   \label{eq: girsanov t 3}
	\end{equation}
        By the same reasoning in Lemma \ref{lem: compensated martingale}, the LHS of \eqref{eq: girsanov t 3} then becomes
        \begin{equation}\begin{aligned}
				\boldsymbol{E}_{\boldsymbol{P}}[\int_{0}^{t} Z_{s-}C_s \lambda_sdN_s] =  \boldsymbol{E}_{\boldsymbol{P}}[\int_{0}^{t}Z_{s-}C_s \lambda_sd{s}]
        			\end{aligned}
                    \label{eq: girsanov t 4}
	\end{equation}
        where we use the fact that \(N(\omega,t) - t\) is a $\boldsymbol{F}^N$-martingale under $\boldsymbol{P}$. \eqref{eq: girsanov t 4} yields the desired result \eqref{eq: girsanov t 3} because any countable set on $[0,t]$ has Lebesgue measure zero. Hence, the restriction $\boldsymbol{Q}_t$ of $\boldsymbol{Q}$ on ${\mathscr{F}_t}$ shares the common compensator with $\tilde{\boldsymbol{Q}}_t$, so \eqref{eq: girsanov transformation} holds for fixed $t\ge 0$. By Lemma \ref{lem: revuz yor}, \eqref{eq: girsanov transformation} holds also for $t =\infty$.
        \end{proof}

\section{Applications to the Cox and Hawkes Processes}
The Cox process, also known as the doubly stochastic Poisson process, is usually introduced as a random measure when some random measure is given. In fact, \cite{prohorov1956convergence} furnished a complete and separable metric on such space of measures, which permits conditioning. Here we adopt a different definition.  Suppose $\underline{\underline{N}} = \{N(\omega,t),\boldsymbol{F}, \boldsymbol{P}\}$ satisfies \textnormal{(P0)-(P5)}. A Cox process $\underline{\underline{N}}$ characterized is by its compensator $A(\omega,t)$ in the following way
    \begin{enumerate}
        \item[(C0)] (adapted) $N(\omega,t),A(\omega,t)$ is adapted to $\boldsymbol F$;
        \item[(C1)] (initial adapted) $\mathscr{F}_{0}$ contains the $\sigma$-algebra generated by all $A(\cdot,t), t \geq 0$;
        \item[(C2)] (conditional Poisson law) the random variable $N(\cdot,t) - N(\cdot,s), 0 \leq s < t$ has the conditional Poisson law for $n =0, 1,2,\dots$,
		\[\boldsymbol{P}(N(\cdot,t) - N(\cdot,s) = n|\mathscr{F}_{s})(\omega) = \frac{(A(\omega,t) - A(\omega,s))^{n}e^{A(\omega,s) - A(\omega,t)}}{n!};\]
		\item[(C3)] (conditional independent increment) the random variables $N(\cdot,t) - N(\cdot,s)$ and $N(\cdot,u) - N(\cdot,v)$ are conditionally independent given $\mathscr{F}_{0}$ for $0 \leq s < t \leq u < v$.
    \end{enumerate}

        From Theorem \ref{thm: main theorem 1}, Theorem \ref{thm: main theorem 2} and the Lebesgue decomposition theorem, we have 

        \begin{theorem}\label{thm: decomposition theorem}
            Any Cox process decomposes into three parts
            \begin{equation}
                N + A + M
                \label{eq: decomposition theorem}
            \end{equation}
            where $N$ is a Cox process admitting an intensity, $A$ is singular with respect to the Lebesgue measure almost surely, and $M$ is a martingale. In particular, if $N$ is independent of $A$, then $N$ is a non-homogeneous Poisson process
        \end{theorem}
        \begin{lemma}\label{thm: non-anticipative Poisson process}
            A Cox process with a natural intensity is a Poisson process. 
        \end{lemma}
        \begin{proof}[Proof of Lemma \ref{thm: non-anticipative Poisson process}]
            If a Cox process $\underline{\underline{N}} = \{N(\omega,t), \boldsymbol{F}, \boldsymbol{P}\}$  has a natural intensity $\lambda(\omega,t)$, it is non-anticipative by Theorem \ref{thm: main theorem 1}. As we see that for each $t\ge 0$, $\lambda(\cdot, t)$ is measurable with respect to $\mathscr{F}_{0} = \{\emptyset,\Omega\}$, the intensity $\lambda(\omega,t)$ is therefore deterministic. This implies that $N(\omega,t)$ is a Poisson process.
        \end{proof}
        
        \begin{lemma}\label{lem: lebesgue decomposition}
            The absolutely continuous and singular components in the Lebesgue decomposition of an adapted, right-continuous, and non-decreasing process are both adapted.
        \end{lemma}
        \begin{proof}[Proof of Lemma \ref{lem: lebesgue decomposition}]
        Let $D[0,\infty)$ be the space of all c\`adl\`ag functions on $[0,\infty)$ equipped with the metrizable topology in the sense of \cite{skorokhod1956limit} (See Section 12 in \cite{billingsley2013convergence}). Define $F(\omega) \subseteq D[0,\infty) \times D[0,\infty)$ to be all pairs $(B,C)$ in the Lebesgue decomposition $A = B+C$ of an adapted, right-continuous, and non-decreasing process $A$, where $B$ is absolutely continuous and $C$ is singular with respect to the Lebesgue measure. For fixed $\omega \in \boldsymbol \varOmega$, $F(\omega)$ is non-empty. By the measurable section theorem (Theorem 81 of Chapter III in \cite{Dellacherie1980probabilities}), we can choose a measurable process $(B,C):\boldsymbol \varOmega \to D[0,\infty) \times D[0,\infty)$ such that the Lebesgue decomposition holds. 
            \begin{figure}[htbp]
		\centering
		\[\begin{tikzcd}
			\boldsymbol \varOmega \arrow{r}{(B,C)} \arrow[swap]{dr}{B_t } & D[0,\infty) \times D[0,\infty)
			\arrow{d}{\pi^1_t} \\
			& {[0,\infty)}
		\end{tikzcd}\]
	\end{figure}
    The fact that $B$ is adapted follows from the commutative diagram, where $\pi_t^1: (B_{\cdot},C_{\cdot}) \mapsto B_t$ is the measurable projection from $D[0,\infty) \times D[0,\infty)$ to ${[0,\infty)}$. Likewise, we can use another measurable projection $\pi^2_t: (B_{\cdot},C_{\cdot}) \mapsto C_t$ from $D[0,\infty) \times D[0,\infty)$ to ${[0,\infty)}$ and a similar commutative relation to conclude that $C$ is also adapted.
        \end{proof}
        
        \begin{proof}[Proof of Theorem \ref{thm: decomposition theorem}] Let $\underline{\underline{(N, A)}} = \{N(\omega,t), A(\omega,t), \boldsymbol{F}, \boldsymbol{P}\}$ be a Cox process. If $\boldsymbol {F}^{N}$ is the natural filtration generated by $N(\cdot,s), 0 \le s \le t$ and $A^{N}(\omega,t)$ is the predictable projection of $A(\omega,t)$ onto $\boldsymbol{F}^N$ (See the definition in Remark 2.), then $N(\omega,t) - A^{N}(\omega,t)$ is an $\boldsymbol{F}^N$-martingale. By Lemma \ref{lem: lebesgue decomposition}, there exists an $\boldsymbol F^N$-predictable set $E \subseteq \boldsymbol \varOmega\times [0,\infty)$ on which $N(\omega,t)$ decomposes as
        \[N^s(\omega,t) = \int_0^t \mathbf{1}_E(\omega,s) N(\omega,ds), \quad N^c(\omega,t) = N(\omega,t)  - N^s(\omega,t).\]
        Moreover, the point process $N^s(\omega,t)$ has an $\boldsymbol F^N$-predictable compensator $A^{N,s}(\omega,t)$ according to the Lebesgue singular component. Applying the Doob-Meyer decomposition theorem to $N^{s}(\omega,t)$, this yields the decomposition \eqref{eq: decomposition theorem}. In addition, if $N^c(\omega,t)$ is independent of $A^s(\omega,t)$, then the distribution of $N^c(\omega,t)$ is independent of $E$ since $E$ is the support of $A^s(\omega,t)$. This implies that the $\boldsymbol F^N$-predictable set $E$ decomposes as $E =  \boldsymbol \varOmega \times E'$, where $E'$ is a subset of $[0,\infty)$. Furthermore, $E$ is $\boldsymbol F^{N^c}$-predictable, and Lemma \ref{thm: non-anticipative Poisson process} applies to $N^{c}(\omega,t)$. That is, $N^c(\omega,t)$ is a non-homogeneous Poisson process.
        \end{proof}

        
		Suppose $\underline{\underline{N}} = \{N(\omega,t),\boldsymbol{F}, \boldsymbol{P}\}$ satisfies \textnormal{(P0)-(P5)}. A Hawkes process $N(\omega,t)$ has arrivals at times $0(=T_{0}(\omega))<T_1(\omega)<T_{2}(\omega)<\cdots$ such that its intensity satisfies
		\begin{equation}\label{eq: Hawkes}
			\lambda(\omega,t)= \mu (t)+\sum _{{n}:T_n(\omega)< t}\phi (t-T_{n}(\omega)).
		\end{equation} 
		where $\mu$ is a non-negative and locally integrable function on $[0,\infty)$, called the underlying baseline intensity, and $\phi$ is a non-negative and continuous function on $[0,\infty)$ satisfying 
        \[m = \int_{0}^{\infty} \phi(t)dt < \infty\]
        and $\phi(t) = 0$ for $t < 0$, namely the self-exciting function. 
        In general, this process is non-anticipative and may have jumps dependent on its history. \cite{hawkes1974} proved the existence of such processes in the mean-field limit sense. By Theorem \ref{thm: main theorem 1}, we could study their connections to the Poisson process. 
	
        \begin{lemma}\label{lem: constant point process}
            A point process $N(\omega,t)$ is constant almost surely on the intervals 
            \[(s,t), \quad [s,t),\quad  (s,t],\quad [s,t], \quad s<t\]
            if and only if so is its compensator $A(\omega,t)$.
        \end{lemma}
        \begin{proof}[Proof of Lemma \ref{lem: constant point process}]
           Let $\bar N(\omega,t) = N(\omega,t) - N(\omega,s)$. If $N(\omega,t)$ is constant almost surely, $\bar N(\omega,t)$ is a point process $\equiv$ 0. By the Doob-Meyer decomposition theorem, the compensator $\equiv$ 0. If the compensator $A(\omega,t)$ is constant almost surely, then $\bar N(\omega,t)$ is a point process with compensator $\bar A(\omega,t) \equiv 0$. By the statement (1) in Lemma \ref{thm: watanabe-bremaud}, $\bar N(\omega,t) \equiv 0$.
        \end{proof}
        
        \begin{theorem}\label{thm: Hawkes Poisson process}
            A Cox-Hawkes process has a zero self-exciting function. 
        \end{theorem}
                        \begin{figure}[!ht]
        \begin{minipage}[htbp]{.49\linewidth}
            \centering
            \subfloat[][$\phi = 0$]{\label{Rplot1}
            \includegraphics[width=1\linewidth]{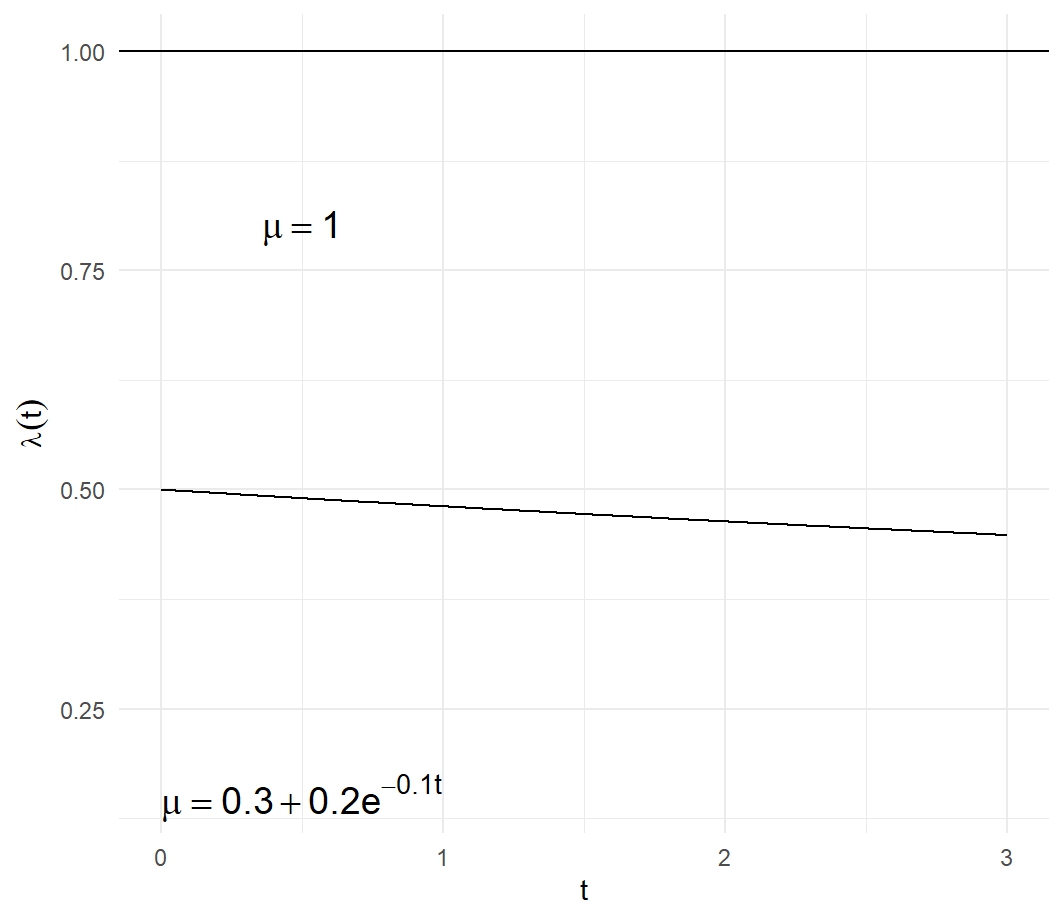}}
        \end{minipage} 
        \hfill    
        \begin{minipage}[htbp]{.49\linewidth}
            \centering
            \subfloat[][$\phi\not\equiv 0$]{\label{Rplot2}
            \includegraphics[width=1\linewidth]{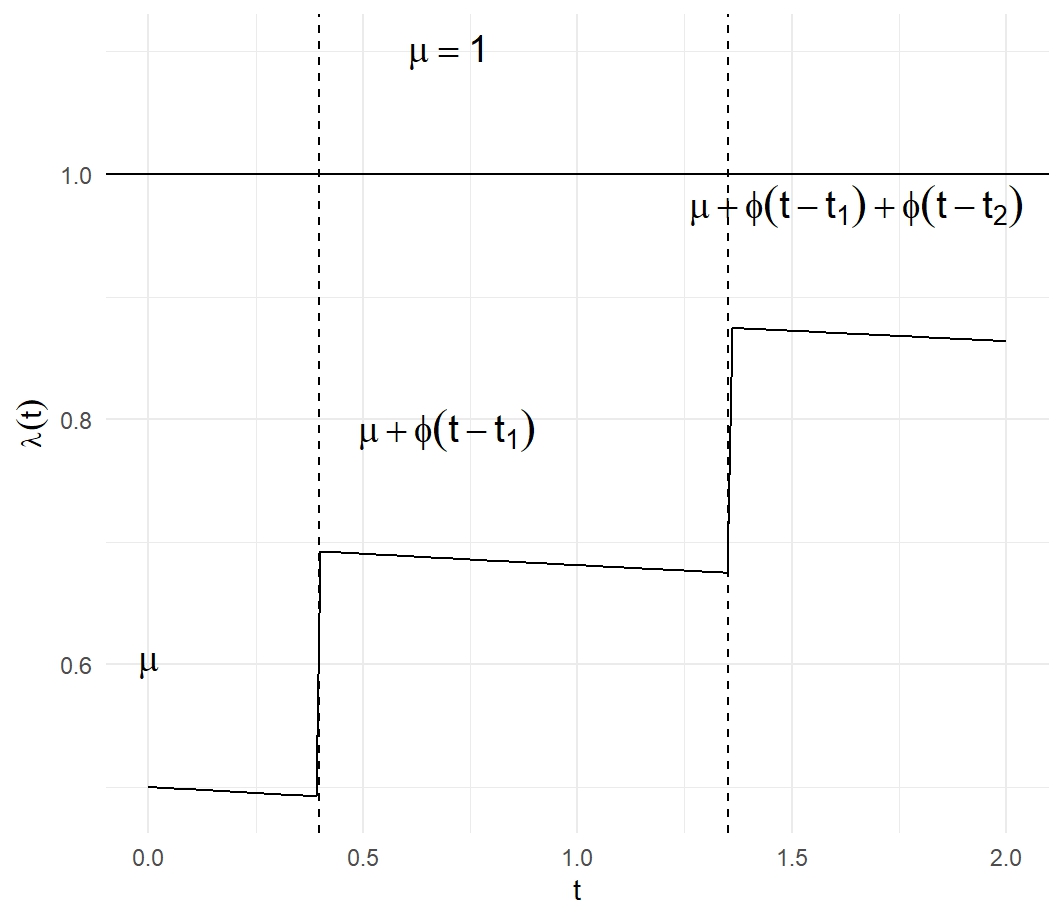}}
        \end{minipage}
        \caption{The intensity $\lambda_{t}$ over time $t$ in a path of Hawkes process under two scenarios (a) $\mu(t) = 0.3 + 0.2 e^{-0.1t}$ and (b) $\mu(t) = 0.1 + 0.2 e^{-0.1t}$ and $\phi(t) = 0.2$.}
        \label{fig:Genelecs}
    \end{figure}

    \begin{proof}[Proof of Theorem \ref{thm: Hawkes Poisson process}]
            Suppose $\underline{\underline{N}} = \{N(\omega,t),\boldsymbol{F}, \boldsymbol{P}\}$ is a Hawkes process. We will show that $\phi \equiv 0$.  
            If there exists $t_{0}\geq 0$ such that $\phi(t_{0}) \neq 0$, there exists $\delta >0$ such that $\phi > 0$ on $[t_{0}, t_{0} + \delta]$ since $\phi$ is continuous on $[0,\infty)$. The point process $N(\omega,t)$ jumps on the interval $(0,t_{0}+\delta]$ with a positive probability, since if not, $A(\omega,t+\delta)>0$ and it contradicts Lemma \ref{lem: constant point process}. This is impossible, because $N(\omega,t)$ will jump at time $t^{\prime}$ on the interval $(0,t_{0}+\delta]$ while $\phi(t - t^{\prime})\neq 0$ on $[t^{\prime}+t_{0}, t^{\prime} + t_{0} + \delta] \subseteq (t_{0}, 2(t_{0} + \delta)]$ and $\lambda(\omega,t)$ will not be deterministic on $(t_{0}, 2(t_{0} + \delta)]$. 
        \end{proof}


    \bibliographystyle{abbrvnat}
    \bibliography{bibtex} 

\end{document}